\documentclass[a4paper,12pt,reqno]{amsart}
\usepackage{amssymb}
\usepackage[dvips]{graphicx}
\usepackage{hyperref}

\usepackage[notcite]{showkeys} 

\newtheorem{cor}{Corollary}[section]
\newtheorem{nncor}{Corollary}

\input{parab.def}
\input{comments.tex}
\usepackage[notcite]{showkeys} 

\title{Loop space of a K\"ahler manifold}

\author{Anakkar M.}

\begin{document}
	
	\maketitle
	
	\newsect[ABS]{Abstract}
We prove that the loop space of a Kähler manifold inherits a Kähler structure. Then we prove that equipped with this natural metric the loop space is complete and unboudned. Additionally, we show that a geodesic on the loop space can be constructed by piecing together geodesics from each individual leaf.

\newsect[INT]{Introduction}
\subsection{K\"ahlarity of the metric} 
Let $X$ be a $n$-dimensional complex Kähler manifold. Denote by $h$ the K\"ahler metric and let
$\omega$ be its associated Kähler form, i.e $\omega = -\mathrm{Im}(h)$. We consider the Sobolev 
loop space $LX = W^{k,2}(\mathbb{S}^1,X)$. One can define a Hermitian scalar product $H$ on $LX$ by
\begin{equation}
	\eqqno(herm-m1)
	H_g(\xi,\eta) = \int_{\mathbb{S}^1} h_{g(s)}\Big(\xi(s), \eta(s)\Big) \mathrm{d}s.
\end{equation}
Here $g \in LX$ and $\xi,\eta \in T_g(LX)= W^{k,2}(\mathbb{S}^1,g^*TX)$.
The associated with $H$ K\"ahler form $\Omega = -\mathrm{Im}(H)$ writes as 
\begin{equation}
	\eqqno(kahler-f1)
	\Omega_g(\xi,\eta) = \int_{\mathbb{S}^1} \omega_{g(s)}(\xi(s), \eta(s)) \mathrm{d}s.
\end{equation}
Our goal in this paper is to prove the following.
\begin{nnthm}
	Let $\xi$, $\eta$ and $\nu$ three vector fields on $LX$. Then 
\[
d\Omega(\xi, \eta,\nu) = \int_{\mathbb{S}^1}(d\omega)(\xi(s), \eta(s), \nu(s))ds \, .
\]
\end{nnthm}
\begin{nncor} If $(X,h)$ is K\"ahler then 
	$(LX, H)$ is K\"ahler as well.  
\end{nncor} 

An analogous statement  is known in the following two cases. 
\begin{itemize}
	\item[1.]Consider a compact lie group $G$ with an invariant inner product $\langle ., .\rangle$ on the lie algebra $\mathfrak{g}$. Define a form on the smooth loop space $LG = \mathcal{C}^{\infty}(\ss^1, G)$ as follows.
	\begin{equation} \eqqno(SKah)
		\Omega(\xi, \eta) = \frac{1}{2\pi}\int_0^{2\pi} \langle \xi(e^{i\theta}) , \frac{d\eta}{d\theta}(e^{i\theta}) \rangle d \theta.
	\end{equation}
	Then it is proved in \cite{S} that $LG$ equipped with this $\Omega$ is K\"ahler.
	\item[2.] In the Brylinski's book \cite{B} a K\"ahler form is constructed on an oriented three-manifold. More precisely, by using a volume form $\nu$ on the three-manifold $M$ one constructs the K\"ahler form $\Omega$ along the loop $\gamma \in LM$ by
	\begin{equation} \eqqno(BKah)
		\Omega_\gamma(\xi, \eta) = \int_{0}^1 \nu \Big( \frac{d\gamma}{dx}(x), \xi(x), \eta(x) \Big) dx
	\end{equation}
\end{itemize}
Notice that the definition \eqqref(SKah) and \eqqref(BKah) are different from \eqqref(kahler-f1).

\subsection{$L\pp^1$ versus $\pp(l^2)$} Next we study two model of infinite dimensional Hilbert complex manifold : $L\pp^1$ and $\pp(l^2)$. On the one hand, the loop space $L\pp^1$ equipped with the metric defined in \eqqref(herm-m1) from the Fubini-Study metric on $\pp^1$ gives us that $L\pp^1$ is unbounded and complete. On the other hand, $\pp(l^2)$ equipped with the Fubini-Study metric is complete and bounded. 

\subsection{Geodesics in loop spaces} On the loop space $LX$ the natural hermitian metric allow us to define a Levi-Civita connection. We define another connection $\tilde \nabla$ on $LX$ and  
 we prove that this is a connection that is compatible with the metric i.e it is a Levi-Civita connection. By uniqueness of the Levi-Civita connection one can express the geodesic equation by the constructed connection : 
For a path $\gamma$ on $LX$ the geodesic equation is given by $\tilde{\nabla}_{\dot\gamma} \dot\gamma = 0$. Thus we have the following. 

\begin{thm}
	Let $f$ and $g$ two loops on $LX$ and let $\gamma$ be a path joining them, i.e.
	\[
	\forall s \in \ss^1, \quad \gamma_s(0)= f(s) \text{ and } \gamma_s(1)=g(s).
	\]
	Then $\gamma$ is a geodesic in $LX$ if and only if for all $s \in \ss^1$ $\gamma_s$ is a geodesic on $X$ joining $f(s)$ and $g(s)$. 
\end{thm}
\newsect[KAH]{K\"alarity}

\smallskip As it is well known, see \cite{L}, for vector fields $\xi, \eta, 
\nu: LX \to T(LX)$ one has: 
\[
d\Omega(\xi, \eta, \nu) = \xi(\Omega(\eta, \nu))- \eta(\Omega(\xi, \nu)) + 
\nu( \Omega(\xi, \eta)) -
\]
\begin{equation}
	\eqqno(diff-1)
	\qquad\qquad - \Omega([\xi,\eta],\nu) - \Omega([\eta, \nu],\xi) + \Omega([\xi,\nu], \eta).
\end{equation}
First we compute the term $\nu(\Omega(\xi, \eta))(g)$ on the right hand side of \eqqref(diff-1).
Let us notice the following, see \cite{L}. For a function $F$ on $LX$ and $g \in LX$ denote as 
$d_gF$ the De Rham differential of $F$ at $g$. If $\nu$ is a tangent vector to $LX$ and 
$\nu_g \in T_g(LX)= W^{k,2}(\ss^1,g^*TX)$ its value at $g$ we have  
\begin{equation}
	\eqqno(diff-2)
	\nu (F)(g) = (d_gF)[\nu_g],
\end{equation}
where $\nu(F)$ is the derivative of $F$ along $\nu$ and $\nu(F)(g)$ its value at $g$.  
Therefore we have that
\begin{equation}
	\eqqno(diff-3)
	\nu(\Omega(\xi, \eta))(g) = d_g\Big(\Omega(\xi, \eta)\Big) [\nu_g].
\end{equation}

\begin{lem} Under the condition as above one has
\label{Dint=intD}	\begin{equation}
		\eqqno(diff-4)
		d_g\Big(\Omega(\xi, \eta)\Big) [\nu_g]= \int_{\ss^1} \left[ d_g\left( \omega_{g} 
		\big( \xi_g, \eta_g \big) \right) \right] [\nu_g] (s) \, ds.
	\end{equation}

\end{lem}
\proof Consider the linear map $I: W^{k,2}(\ss^1, \cc) \to \cc$ defined as  
\begin{equation}
	I(f) = \int_{\ss^1}f(s)ds \,  .
\end{equation}
Its differential at $f$ in the direction $\nu\in W^{k,2}(\ss^1,\cc)$ is 
\begin{equation}
	\eqqno(diff-i)
	d_fI[\nu ] = \frac{dI\left(f + t \nu \right)}{dt}\Big\vert_{t=0}
	= \frac{d}{dt}\left[ \int_{\ss^1} f(s)+t\nu(s)ds \right]_{t=0} = 
	\int_{\ss^1} \nu(s)ds.
\end{equation}
Then for $g \in LX$ and $\xi, \eta ,\nu \in TLX$ we have 
\[
 d_g\left( \Omega(\xi, \eta) \right)[\nu_g] =
d_g\left[\int_{\ss^1}\omega_{g(s)}(\xi_g (s),\eta_g (s))ds\right][\nu_g] = d_g\Big( I(\omega_g(\xi_g,\eta_g)) \Big)[\nu_g]
\]
\begin{equation}  \eqqno(diff-i2)
	 = d_{\omega_g(\xi_g,\eta_g)}I \Big(
	d_g(\omega_g(\xi_g,\eta_g))[\nu_g] \Big) =
	\int_{\mathbb{S}^1} \left[ d_g(	\omega_{g}(\xi_g,\eta_g)) \right][\nu_g] (s)ds 
\end{equation}
by \eqqref(diff-i). Here $g \mapsto \omega_g(\xi_g,\eta_g)$ we consider as a map 
$W^{k,2}(\ss^1,X) \to W^{k,2}(\ss^1, \cc)$. Its De Rham differential 
$d_g(\omega_g(\xi_g,\eta_g))$ is a 1-form with values in $W^{k,2}(\ss^1, \cc)$.
\smallskip \qed

\begin{lem} Under the condition above one has
\label{evD=Dev}	\begin{equation}
	\forall s \in \ss^1, \quad	d_g\left( \omega_{g} \big( \xi_g, \eta_g \big) \right) [\nu_g](s)= d_g\left( \omega_{g(s)} 
		\big( \xi_g(s), \eta_g(s) \big) \right) [\nu_g].
	\end{equation}

\end{lem}
\proof 
Let $s \in \ss^1$. We consider the map $ev_s: W^{k,2}(\ss^1,\cc) \to \cc $ defined as $ev_s(g) = g(s)$. 
Its differential  
$
d_gev_s:W^{k,2}(\ss^1,\cc) \to  \cc
$
at $g$ taken in the direction $\xi \in T_g(LX)= W^{k,2}(\ss^1, g^*TX)$ is
$  
d_gev_s[\xi]=  \xi(s) $. In fact, we can write locally 
\[
d_gev_s[\xi] = \frac{d \left( ev_s( g+t\xi) \right) }{dt} \Big\vert_{t=0}  = \frac{d \left( g(s) + t\xi(s) \right) }{dt} 
\Big\vert_{t=0} = \xi(s).
\]
Let define $\Theta : LX \to W^{k,2}(\ss^1, \cc)$ by $\Theta(g) = \omega_g(\xi_g, \eta_g)$ and $K_s= ev_s\circ \Theta$, i.e 
\[
\forall g \in LX, \quad K_s(g)= \omega_{g(s)}(\xi_g(s), \eta_g(s)).
\]
 Then we have 
\[
\left[ d_g(	\omega_{g}(\xi_g,\eta_g)) \right][\nu_g] (s) = d_{\Theta(g)}ev_s
\left[ d_g\Theta \right][\nu_g]= d_g\left( ev_s \circ \Theta \right)[\nu_g] =
\]
\begin{equation}
=d_gK_s[\nu_g] = d_g \left(\omega_{g(s)}(\xi_g(s), \eta_g(s))\right) [\nu_g].
\end{equation}
\smallskip \qed
\noindent By combining lemma \ref{Dint=intD} and Lemma \ref{evD=Dev} we obtain the following equality. 
\begin{equation} \eqqno(f-diff-4)
\nu(\Omega(\xi, \eta))(g) = \int_{s \in \ss^1}	d_g\left( \omega_{g(s)} 
\big( \xi_g(s), \eta_g(s) \big) \right) [\nu_g] \, ds.
\end{equation}
For the next computation we need the following technical Lemma.

\begin{lem} For $s\in \ss^1$ we consider the map $T_s: LX \to X$ defined as $T_s(g) = g(s)$. Then 
\[
\forall \xi, \nu \in T_g(LX), \quad d_gT_s[\xi]=  \xi(s) \quad \text{ and } \quad  d^2_gT_s[\xi, \nu]=0
\]	
\end{lem}
\proof
Let $s \in \ss^1$. We consider the map $T_s: LX \to X$ defined as $T_s(g) = g(s)$. 
Its differential $dT_s(g):T_g(LX) \to T_{g(s)}X$
at $g$ taken in the direction $\xi \in T_g(LX)= W^{k,2}(\ss^1, g^*TX)$ is
\[   
d_gT_s[\xi]=  \xi(s) \in  T_{g(s)}X.
\]
Indeed, locally we can write
\[
d_gT_s[\xi] = \frac{d T_s( g+t\xi) }{dt} \Big\vert_{t=0}  = \frac{d \left( g(s) + t\xi(s) \right) }{dt} 
\Big\vert_{t=0} = \xi(s) \in T_{g(s)}X.
\]
So we have   
\[ d^2T_s [\xi, \nu] = \frac{d\left(d_{g +t\nu}T_s[\xi]\right)}{dt}\Big\vert_{t=0} = 
\frac{d \xi(s)}{dt}\Big\vert_{t=0}= 0 . \]
\smallskip \qed

\smallskip Take a chart $(\phi, U)$ on $X$ containing $g(s)$. In the image $\phi(U) = \tilde{U} \subset \cc^n$,  
$\tilde \omega = (\phi^{-1})^*\omega$ is given by 
\[ 
\tilde \omega_z = \frac{i}{2}\sum_{i,j=1}^n h_{ij}(z) dz_i \wedge d\overline{z}_j,  
\]
where $(h_{ij})$ is a Hermitian matrix. Then we have the following expression for the term in the integral \eqqref(f-diff-4). 
\[
	\omega_{g(s)}(\xi_g(s),\eta_g(s)) = \tilde{\omega}_{\phi(g(s))}(d_{g(s)}\phi [\xi_g(s)], d_{g(s)}\phi [\eta_g(s)])
\]
\begin{equation}
\sum_{i,j=1}^nh_{ij}(\phi(g(s)))\, dz_i\wedge d\bar{z}_j \left(  d_{g(s)}\phi [\xi_g(s)], d_{g(s)}\phi [\eta_g(s)] 
\right)
\end{equation}

\begin{lem} One has the following expression for the term under the integral in \eqqref(f-diff-4)
	\[ 
	d_g \Big(\omega_{g(s)}(\xi_g(s),\eta_g(s))\Big)[\nu_g] = A+B-C.
	\] 
	Here
	
	\smallskip\noindent 
	\begin{equation}
		A=\frac{i}{2} \sum_{i,j=1}^n \sum_{k=1}^n \frac{\partial h_{ij}}{\partial z_k}\nu_k\xi_i\overline{\eta_j} +
		\frac{\partial h_{ij}}{\partial \overline{z}_k}\overline{\nu_k}\xi_i\overline{\eta_j} -
		\frac{\partial h_{ij}}{\partial z_k}\nu_k\eta_i\overline{\xi_j} -\frac{\partial h_{ij}}{\partial 
			\overline{z}_k}\overline{\nu_k}\eta_i\overline{\xi_j}; 
	\end{equation}
	As well as 
	\noindent 
	\[ B = \frac{i}{2} \sum_{i,j=1}^n h_{ij}d^2\phi_i\big( \xi,\nu \big) \overline{\eta_j} 
	+ h_{ij}d\xi_i[\nu]\overline{\eta_j} + h_{ij}\xi_i \overline{d^2\phi_j\big( \eta,\nu \big)} 
	+ h_{ij}\xi_i \overline{d\eta_j[\nu]};
	\]

	\noindent 
	\[
	C = \frac{i}{2} \sum_{i,j=1}^n h_{ij}\overline{d^2\phi_j\big( \xi,\nu \big)}\eta_i + h_{ij}\overline{d\xi_j[\nu]}\eta_i + h_{ij}\overline{\xi_j} d^2\phi_i\big( \eta,\nu \big) + h_{ij}\overline{\xi_j} d\eta_i[\nu]
	\]
	where $\zeta_l$ stands for $dz_ld_{g(s)}\phi[\zeta_g(s)]$ with $ l\in \{ i,j,k \}$ and $\zeta \in \{ \xi, \eta, \nu \}$ for simplificity, i.e. 
	\[
	\begin{array}{ccc}
\xi_i = dz_i d_{g(s)}\phi[\xi_g(s)] & \eta_i = dz_i d_{g(s)}\phi[\eta_g(s)] & \nu_i = dz_i d_{g(s)}\phi[\nu_g(s)] \\
\xi_j = dz_j d_{g(s)}\phi[\xi_g(s)] & \eta_j = dz_j d_{g(s)}\phi[\eta_g(s)] & \nu_j = dz_j d_{g(s)}\phi[\nu_g(s)]
\\
\xi_k = dz_k d_{g(s)}\phi[\xi_g(s)] & \eta_k = dz_k d_{g(s)}\phi[\eta_g(s)] & \nu_k = dz_k d_{g(s)}\phi[\nu_g(s)].
	\end{array}
	\] 
	As well  $d^2\phi_l(\zeta, \zeta')$ stands for $dz_ld^2_{g(s)}\phi(\zeta_g(s), \zeta'_g(s))$  and $ d\zeta_l[\zeta']$ stands for $dz_l d_{g(s)}\phi d_g\zeta[\zeta'](s)$ 
	with $\zeta, \zeta' \in \{\xi, \eta, \nu\}$ and $l \in \{i,j,k\}$;
\end{lem}
\proof Write $\omega_{g(s)}(\xi_g(s),\eta_g(s))=	  (\phi^*\tilde\omega)_{g(s)}(\xi_g(s),\eta_g(s)) =$
\[
= \tilde \omega_{\phi(T_s(g))}\big(d_{T_s(g)}\phi [d_gT_s[\xi_g]], d_{T_s(g)}\phi \left[d_gT_s [\eta_g] \right]\big) \]
\[ = \frac{i}{2} \sum_{i,j=1}^n h_{ij}\big( \phi(T_s(g))\big) dz_i\wedge d\overline{z}_j\Big(d_{T_s(g)}\phi [d_gT_s[\xi_g]], d_{T_s(g)}\phi \left[d_gT_s [\eta_g] \right] \Big) 
\]
 \begin{equation} \eqqno(firstsum)
	= \frac{i}{2} \sum_{i,j=1}^n h_{ij}\big( \phi(T_s(g))\big) dz_i \Big( d_{T_s(g)}\phi [d_gT_s[\xi_g]] \Big) d\overline{z}_j\Big( d_{T_s(g)}\phi [d_gT_s [\eta_g]]\Big) -
\end{equation}

\begin{equation} \eqqno(secondsum)
	- \frac{i}{2} \sum_{i,j=1}^n h_{ij}\big( \phi(T_s(g))\big)  d\overline{z}_j\Big( d_{T_s(g)}\phi [d_gT_s[\xi_g]]\Big) dz_i\Big(  d_{T_s(g)}\phi [d_gT_s [\eta_g]]  \Big)
\end{equation} 
For the first sum \eqqref(firstsum) we have to diffrentiate with respect to $g$ a product of three functions (the second sum \eqqref(secondsum) is computed similarly):
\begin{enumerate}
	\item $D_{ij} : g \mapsto h_{ij}\big( \phi(T_s(g))\big)$.
	\item $E_{i} : g\mapsto  dz_i \Big( d_{T_s(g)}\phi [d_gT_s[\xi_g]] \Big)$.
	\item $F_{j} : g \mapsto d\overline{z}_j\Big( d_{T_s(g)}\phi [d_gT_s [\eta_g]]\Big)$.
\end{enumerate}	 
Using the chain rule on the differentiation by $g$ on the direction $\zeta \in T_g(LX)$ we obtain 
\begin{enumerate}
	\item $	d_g  h_{ij}\big( \phi(T_s(g))\big) [\zeta] = $
	\begin{eqnarray*}
		&=& \sum_{k=1}^n \frac{\partial h_{ij}}{\partial z_k}\big( \phi(T_s(g))\big)dz_k \left[ d_{T_s(g)}\phi\left[d_gT_s [\zeta ]\right] \right]  +\sum_{k=1}^n \frac{\partial h_{ij}}{\partial \overline{z}_k}\big( \phi(T_s(g))\big) d\overline{z}_k \left[d_{T_s(g)}\phi \left[d_gT_s [\zeta ]\right] \right] \\
		& = & \sum_{k=1}^n \frac{\partial h_{ij}}{\partial z_k}\big( \phi(g(s))\big) dz_k \left[d_{g(s)}\phi [\zeta(s)]\right]  +\sum_{k=1}^n \frac{\partial h_{ij}}{\partial \overline{z}_k}\big( \phi(g(s))\big) d\overline{z}_k \left[d_{g(s)}\phi [\zeta(s)] \right] \\
		& = & \sum_{k=1}^n \frac{\partial h_{ij}}{\partial z_k}d_{g(s)}\phi_k [\zeta(s)]  +\sum_{k=1}^n \frac{\partial h_{ij}}{\partial \overline{z}_k}\overline{ d_{g(s)}\phi_k [\zeta(s)]} \\
		&=& \sum_{k=1}^n \frac{\partial h_{ij}}{\partial z_k}\zeta_k +\sum_{k=1}^n \frac{\partial h_{ij}}{\partial \overline{z}_k}\overline{\zeta}_k.
	\end{eqnarray*}	

	\item 
	\begin{eqnarray*}
		d_g \left( dz_i \Big( (d\phi)_{T_s(g)}(dT_s)_g\xi_g\Big) \right) [\zeta] &=& dz_i \left( d_{T_s(g)}^2\phi\big( d_gT_s [\xi_g], d_gT_s [\zeta] \big) + 	d_{T_s(g)}\phi \left[d_gT_s \big[d_g\xi [\zeta]\big]\right] \right)  \\
		&=& dz_i \left( d_{g(s)}^2\phi \big( \xi_g(s),\zeta(s) \big) + 	d_{g(s)}\phi \big[d_g\xi [\zeta (s)]\big]\right) \\
		&=& d_{g(s)}^2\phi_i\big( \xi_g(s),\zeta(s) \big) + 	d_{g(s)}\phi_i \big[d_g\xi [\zeta (s)] \big] \\
		&=& d^2\phi_i\big( \xi,\zeta \big) + d\xi_i [\zeta] 
	\end{eqnarray*}
	\item By the same way as 2, 
	\begin{eqnarray*}
		d_g d\bar{z}_j\Big( d_{T_s(g)}\phi [d_gT_s [\eta_g]]\Big) [\zeta] &=& d\bar{z}_j \left( d_{T_s(g)}^2\phi\big( d_gT_s [\eta_g], d_gT_s [\zeta] \big) + 	d_{T_s(g)}\phi \big[d_gT_s [d_g\eta [\zeta]]\big] \right)  \\
		&=& d\bar{z}_j \left( d_{g(s)}^2\phi \big( \eta_g(s),\zeta(s) \big) + 	d_{g(s)}\phi \big[ d_g\eta [\zeta] \big] (s)  \right) \\
		&=& \overline{ d_{g(s)}^2\phi_j \big( \eta_g(s),\zeta(s) \big)} + 	\overline{ d_{g(s)}\phi_j [d_g\eta [\zeta] (s)]} \\
		&=& \overline{ d^2\phi_j \big( \eta,\zeta \big)} + 	\overline{ d\eta_j [\zeta]}
	\end{eqnarray*}
Then we obtain 
\[ d_g\left( h_{ij}\big( \phi(T_s(g))\big) dz_i \Big( d_{T_s(g)}\phi [d_gT_s[\xi_g]] \Big) d\overline{z}_j\Big( d_{T_s(g)}\phi [d_gT_s [\eta_g]]\Big) \right)= \]
\[ =
d_g  h_{ij}\big( \phi(T_s(g))\big) [\zeta] dz_i \Big( d_{T_s(g)}\phi [d_gT_s[\xi_g]] \Big)d\overline{z}_j\Big( d_{T_s(g)}\phi [d_gT_s [\eta_g]]\Big)
+
\]
\[+ h_{ij}\big( \phi(T_s(g))\big)  \, \, d_g\left(dz_i \Big( d_{T_s(g)}\phi [d_gT_s[\xi_g]] \Big) \right)[\zeta] \, d\overline{z}_j\Big( d_{T_s(g)}\phi [d_gT_s [\eta_g]]\Big)+
\]
\[
+ h_{ij}\big( \phi(T_s(g))\big) dz_i \Big( d_{T_s(g)}\phi [d_gT_s[\xi_g]] \Big) d_g\left(d\overline{z}_j\Big( d_{T_s(g)}\phi [d_gT_s [\eta_g]]\Big) \right)
\]
\end{enumerate}
Then we compute 
\begin{equation} \eqqno(dgDEF)
	d_g(D_{ij}E_iF_j) = d_gD_{ij}E_iF_j + D_{ij}d_gE_i F_j + D_{ij}E_id_gF_j.
\end{equation}
Summing all the terms we obtain our Lemma. 
\smallskip\qed

According to \eqqref(diff-1) and \eqqref(diff-4) we have  
\begin{eqnarray*}
	&\textbf{d}\Omega(\xi, \eta, \nu) =  
	\xi(\Omega(\eta, \nu))- \eta(\Omega(\xi, \nu)) + \nu( \Omega(\xi, \eta)) - \Omega([\xi,\eta],\nu) - \Omega([\eta, \nu],\xi) + \Omega([\xi,\nu], \eta) \\
	&= \displaystyle \int_{s \in \ss^1} d (\omega_{g(s)}(\eta_g(s),\nu_g(s)))\xi_g- d (\omega_{g(s)}(\xi_g(s),\nu_g(s)))\eta_g + d (\omega_{g(s)}(\xi_g(s),\eta_g(s)))\nu_g -
\end{eqnarray*}
\begin{equation}
	-\omega_{g(s)}([\xi, \eta](s), \nu(s)) - \omega_{g(s)}([\eta, \nu](s), \xi(s)) + \omega_{g(s)}([\xi, \nu](s), \eta(s))  ds  \quad 
	\eqqno(diff-5)
\end{equation}

\begin{thm}
	For every $s \in \ss^1$ the expression under the integral \eqqref(diff-5) is equal to $(\mathrm{d}\omega)_{g(s)}(\xi_g(s), \eta_g(s),\nu_g(s))$, i.e.
	\[
	\mathrm{d}\Omega(\xi, \eta, \nu) = \int_{s\in \mathbb{S}^1} (\mathrm{d}\omega)_{g(s)}(\xi_g(s),\eta_g(s),\eta_g(s))\, ds.
	\]
\end{thm}
\proof
Computing  every terms of \eqqref(diff-5) as in the previous lemma we obtain (the terms in the corresponding brackets cancel out)
\[
d (\omega_{g(s)}(\xi_g(s),\eta_g(s)))\nu_g = 
\frac{i}{2} \sum_{i,j=1}^n \sum_{k=1}^n \frac{\partial h_{ij}}{\partial z_k}\nu_k\xi_i\overline{\eta_j} +\frac{\partial h_{ij}}{\partial \overline{z}_k}\overline{\nu_k}\xi_i\overline{\eta_j} -\frac{\partial h_{ij}}{\partial z_k}\nu_k\eta_i\overline{\xi_j} -\frac{\partial h_{ij}}{\partial \overline{z}_k}\overline{\nu_k}\eta_i\overline{\xi_j} 
\]
\[
+\frac{i}{2} \sum_{i,j=1}^n \underbrace{h_{ij} d^2\phi_i\big( \xi,\nu \big) \overline{\eta_j}}_{1} + \underbrace{h_{ij}d\xi_i[\nu]\overline{\eta_j}}_{A} +\underbrace{ h_{ij}\xi_i \overline{d^2\phi_j \big( \eta,\nu \big)}}_{2} + \underbrace{h_{ij}\xi_i \overline{d\eta_j[\nu]}}_{B}
\]
\[ 
+\frac{i}{2} \sum_{i,j=1}^n -\underbrace{h_{ij}\overline{d^2\phi_j\big( \xi,\nu \big)}\eta_i}_{3} - \underbrace{h_{ij}\overline{d\xi_j[\nu]}\eta_i}_{H} - \underbrace{h_{ij}\overline{\xi_j} d^2\phi_i\big( \eta,\nu \big)}_{4} - \underbrace{h_{ij}\overline{\xi_j} d\eta_i[\nu]}_{D}
\]

The same as previous, we compute 
\[
 d (\omega_{g(s)}(\eta_g(s),\nu_g(s)))\xi_g = 
\frac{i}{2} \sum_{i,j=1}^n \sum_{k=1}^n \frac{\partial h_{ij}}{\partial z_k}\xi_k\eta_i\overline{\nu_j} +\frac{\partial h_{ij}}{\partial \overline{z}_k}\overline{\xi_k}\eta_i\overline{\nu_j} -\frac{\partial h_{ij}}{\partial z_k}\xi_k\nu_i\overline{\eta_j} -\frac{\partial h_{ij}}{\partial \overline{z}_k}\overline{\xi_k}\nu_i\overline{\eta_j} 
\]
\[
+\frac{i}{2} \sum_{i,j=1}^n \underbrace{h_{ij}d^2\phi_i\big( \eta,\xi \big) \overline{\nu_j}}_{5} + \underbrace{h_{ij}d\eta_i[\xi]\overline{\nu_j}}_{C} + \underbrace{h_{ij}\eta_i \overline{d^2\phi_j\big( \nu,\xi \big)}}_{3} + \underbrace{h_{ij}\eta_i \overline{d\nu_j[\xi]}}_{I}
\]
\[
+\frac{i}{2} \sum_{i,j=1}^n -\underbrace{h_{ij}\overline{d^2\phi_j\big( \eta,\xi \big)}\nu_i}_{6} - \underbrace{h_{ij}\overline{d\eta_j[\xi]}\nu_i}_{J} - \underbrace{h_{ij}\overline{\eta_j} d^2\phi_i\big( \nu,\xi \big)}_{1} - \underbrace{h_{ij}\overline{\eta_j} d\nu_i[\xi]}_{G}
\]
and the third term
\[ 
d (\omega_{g(s)}(\xi_g(s),\nu_g(s)))\eta_g = 
\frac{i}{2} \sum_{i,j=1}^n \sum_{k=1}^n \frac{\partial h_{ij}}{\partial z_k}\eta_k\xi_i\overline{\nu_j} +\frac{\partial h_{ij}}{\partial \overline{z}_k}\overline{\eta_k}\xi_i\overline{\nu_j} -\frac{\partial h_{ij}}{\partial z_k}\eta_k\nu_i\overline{\xi_j} -\frac{\partial h_{ij}}{\partial \overline{z}_k}\overline{\eta_k}\nu_i\overline{\xi_j} 
\]
\[
+\frac{i}{2} \sum_{i,j=1}^n \underbrace{h_{ij}d^2\phi_i\big( \xi,\eta \big) \overline{\nu_j}}_{5} + \underbrace{h_{ij}d\xi_i[\eta]\overline{\nu_j}}_{E} +\underbrace{ h_{ij}\xi_i \overline{d^2\phi_j\big( \nu,\eta \big)}}_{2} + \underbrace{h_{ij}\xi_i \overline{d\nu_j[\eta]}}_{L}
\]
\[
+\frac{i}{2} \sum_{i,j=1}^n -\underbrace{h_{ij}\overline{d^2\phi_j\big( \xi,\eta \big)}\nu_i}_{6} - \underbrace{h_{ij}\overline{d\xi_j[\eta]}\nu_i }_{K}- \underbrace{h_{ij}\overline{\xi_j} d^2\phi_i\big( \nu,\eta \big)}_{4} - \underbrace{h_{ij}\overline{\xi_j} d\nu_i[\eta]}_{F}
\]
The commutator terms 
\begin{eqnarray}
	\omega_{g(s)}([\xi, \eta](s), \nu(s)) =\frac{i}{2}\sum_{i,j=1}^n \underbrace{h_{ij}d\eta_i[\xi]\overline{\nu_j}}_{C} -\underbrace{ h_{ij}d\xi_i[\eta]\overline{\nu_j} }_{E}
	-\underbrace{h_{ij}\overline{d\eta_j[\xi]}\nu_i }_{J}+ \underbrace{h_{ij}\overline{d\xi_j[\eta]}\nu_i}_{K}
	\\
 \omega_{g(s)}([\eta, \nu](s), \xi(s)) = \frac{i}{2}\sum_{i,j=1}^n \underbrace{h_{ij}d\nu_i[\eta]\overline{\xi_j}}_{F} - \underbrace{h_{ij}d\eta_i[\nu]\overline{\xi_j}}_{D} 
 -\underbrace{h_{ij}\overline{d\nu_j[\eta]}\xi_i}_{L} + \underbrace{h_{ij}\overline{d\eta_j[\nu]}\xi_i}_{B}
 \\
	\omega_{g(s)}([\xi, \nu](s), \eta(s)) = \frac{i}{2}\sum_{i,j=1}^n \underbrace{h_{ij}d\nu_i[\xi]\overline{\eta_j}}_{G} - \underbrace{h_{ij}d\xi_i[\nu]\overline{\eta_j}}_{A}
	-\underbrace{h_{ij}\overline{d\nu_j[\xi]}\eta_i}_{I} + \underbrace{h_{ij}\overline{d\xi_j[\nu]}\eta_i}_{H}
\end{eqnarray}
Then computing the sum according to the sign in \eqqref(diff-5) gives us 
$$ \textbf{d}\Omega(\xi, \eta, \nu) = \int_{\mathbb{S}^1}(\textbf{d}\omega)(\xi(s), \eta(s),\nu(s)))ds$$	
\smallskip\qed
\noindent Then we obtain the following corollary.
\begin{cor}
 If $(X,h)$ is K\"ahler then 
$(LX, H)$ is K\"ahler as well.
\end{cor}
\proof 
Since $\mathrm{d}\omega = 0$ then $\mathrm{d}\Omega =0$. 
\smallskip \qed
\newsect[ex-LP1]{$L\pp^1$ versus $\pp(l^2)$ model.}

Let us consider the loop space $L\pp^1 = W^{1,2}(\ss^1, \pp^1)$ on the Riemann sphere. Let $h$ be the Herimitian metric associated to the Fubini-Study metric. We can endow $L\pp^1$ by the Hermitian metric constructed over the Fubini-Study metric defined as \eqqref(herm-m1). 
By completeness of $\pp^1$ the Hilbert manifold $L\pp^1$ is complete. Let us prove that it is not bounded. Consider two loops $f,g \in L\pp^1$ defined by $f:e^{i\theta}\in \ss^1 \mapsto [1:0]$ and $g:e^{i\theta} \mapsto [\cos(n\theta) :  \, \sin(n\theta)]=[1:\tan(n\theta)]$. 
The geodesic joining $f$ and $g$ may not exists since $LX$ is infinite dimensional, see \cite{L2}. Nevertheless we can compute the geodesic distance by taking the infimum over the paths joining $f$ and $g$.  
Let consider such a path $\gamma$, i.e.  
\[
\forall s \in \ss^1, \quad \gamma_s(0)= f(s) \text{ and } \gamma_s(1)=g(s).
\]
The length of the path $L(\gamma)$ is given by 
\begin{eqnarray*}
	L(\gamma)&=& \int_0^1 \sqrt{H_{\gamma(t)}( \dot{\gamma}(t),\dot{\gamma}(t) ) } dt
	= \int_0^1 \sqrt{ \frac{1}{2\pi} \int_{s \in \ss^1}  h_{\gamma_s(t)}(\dot{\gamma}_s(t),\dot{\gamma}_s(t))ds } dt \\
	& \geqslant & \int_0^1  \int_{s \in \ss^1} \sqrt{  h_{\gamma_s(t)}(\dot{\gamma}_s(t),\dot{\gamma}_s(t)) }ds dt
	  \quad (\text{ by Cauchy-Schwarz inequality })\\
 & = &  \int_{s \in \ss^1} \int_0^1  \sqrt{  h_{\gamma_s(t)}(\dot{\gamma}_s(t),\dot{\gamma}_s(t)) }ds dt
 = \int_{s \in \ss^1}l(\gamma_s)ds
\end{eqnarray*}
where $l(\gamma_s)$ stands for the length of the path $\gamma_s$ joining $f(s)$ to $g(s)$ in $\pp^1$. Let denote  by $\mathrm{dist}_{\pp^1}$ the distance on $\pp^1$ given by the Fubini-Study metric $h$. We have that 
\[
\forall s \in \ss^1, \quad l(\gamma_s) \geqslant \mathrm{dist}_{\pp^1}(f(s), g(s)).
\]
Then we can finish the computing. 
\begin{eqnarray*}
	L(\gamma) & \geqslant & \int_{\mathbb{S}^1} \mathrm{dist}_{\pp^1}(f(s), g(s)) ds 
	= 
	\frac{1}{2\pi}  \int_{\theta = 0}^{2\pi} \int_0^1 \frac{\left|\tan(n\theta)\right|}{1+t^2\left( \tan( n \theta) \right)^2}  dt d\theta\\
	& \geqslant & \frac{1}{{2\pi}} 
	\sum_{k=0}^{4n-1}\int_{\frac{k\pi}{2n}}^{\frac{(k+1)\pi}{2n}}
	\int_0^1 \frac{\left|\tan(n\theta)\right|}{1+t^2\left( \tan( n \theta) \right)^2}  dt d\theta \\
	& \geqslant & \frac{1}{2\pi} 
	\sum_{k=0}^{4n-1}\int_{\frac{k\pi}{2n}}^{\frac{(k+1)\pi}{2n}}
	\int_0^1 \frac{(-1)^k\tan(n\theta)}{1+t^2\left( \tan( n \theta) \right)^2}  dt d\theta \\
	& \geqslant & \frac{1}{2\pi} 
	\sum_{k=0}^{4n-1}\int_{\frac{k\pi}{2n}}^{\frac{(k+1)\pi}{2n}}
	(-1)^k\arctan\left( \tan(n \theta)\right) d\theta \\
	& \geqslant & \frac{1}{2\pi} 
	\sum_{k=0}^{4n-1} \frac{\pi}{2} = \frac{2n \pi}{2\pi}=n 
\end{eqnarray*}
By taking the infimum length $L(\gamma)$ over all the paths $\gamma$ in $L\pp^1$ joining $f$ to $g$ we obtain
\begin{equation}
	\mathrm{dist}_{L\pp^1}(f,g) \geqslant n \xrightarrow[n \to +\infty]{} + \infty
\end{equation}
So $L\pp^1$ is not bounded. 

Now let us study the manifold $\pp(l^2)$. Endowed with its Fubini-Study metric, $\pp(l^2)$ is bounded. 
Recall the Fubini-Study form in $\pp(l^2)$
$$\omega=\frac{i}{2} \left( \frac{1}{\norm{w}^2} \sum_{k=0}^\infty dw_k \wedge d\bar w_k  - \frac{1}{\norm{w}^4} \sum_{k=0}^\infty \bar w_kdw_k \wedge  \sum_{k=0}^\infty w_k d\bar w_k  \right) =dd^c\ln(\norm{w}^2)$$
where $w$ are homogeonous coordinates.
The unitary group $\calu(l^2) = \{ L \in \call(l^2) \ | \ LL^*=\text{Id} \}$ acts holomorphically and transitively on $\pp(l^2)$. Morevoer $\calu(l^2)$ preserves the Fubini-Study form. 
Take two points $[p]$ and $[q]$ in $\pp(l^2)$ and consider a chart $\calu$ that contains these two points. One can suppose that $\calu = \calu_0$ and $\norm{p}=1$. 

Consider an orthonormal basis of $\text{Vect}(p,q)^\perp$ denoted by $(v_j)_{j \in [\![ 3, +\infty [\![}$ and use the Gram-Schmidt algorithm on $p$ and $q$.
Then $v_1= p$ and $v_2 = \alpha( q - <q,p> p )$ with $\alpha = \norm{ q - <q,p> p}^{-1}$. Therefore $(v_j)_{j\in \nn^*}$ is an orthonormal Hilbert basis of $l^2$. Let $(e_j)_{j \in \nn^*}$ be the canonical basis of $l^2$ and then we can define $L \in \call(l^2)$ by $L(e_i)=v_j$. $L$ is unitary and $L^{-1}(q)= \gamma e_2 + \delta e_1 $. In $\pp(l^2)$ coordinates 
$L^* \cdot [p] = [L^*p]=[1 : 0: \dots]$ and $L^* \cdot [q] = [L^*q]=[\gamma : \delta: 0 : \dots]=[1 : \frac{\delta}{\gamma}: 0 : \dots]$ if $\gamma \neq 0$. The case $\gamma = 0$  means that $p$ and $q$ are orthogonal. Then we can take a sequence of points $(q_n)_n$ non orthogonal to $p$ and converging to $q$. Then we replace $q$ by $q_n$ in the following and compute the limit at the end. 

Therefore in $\calu_0$ we brought $p$ to the origin and $q$ in Vect$(e_1)$. Consider the path $\gamma$ defined by $\gamma(t)= t$ for $t \in [0,R]$ where $\gamma(R)=q$. On the chart $\calu_0$ one has $\omega_0= dd^c\ln (1 + \sum_{k=1}^\infty |z_k^0|^2) $. 
$$\text{dist}([p],[q])= \int_0^R\sqrt{\omega_{\gamma(t)}(\gamma'(t),\gamma'(t))}dt =
\int_0^R \sqrt{ \frac{1}{1+t^2} - \frac{t^2}{(1+t^2)^2} } dt= \int_0^R \frac{dt}{1+t^2}$$
Then $\text{dist}(p,q)= \arctan(R) \leqslant \frac{\pi}{2}$. So $\pp(l^2)$ is bounded. 

Second, we remark that
$\pp(l^2)$ is complete. Take $([w_n])_n$ a Cauchy sequence in $\pp(l^2)$. Then there exists $n_0$ such that for every $n,m \geqslant n_0$ the distance $\text{dist}([w_n],[w_m])< \frac{1}{4}$. 
Therefore for every $n\geqslant n_0$ the point $[w_n]$ is in the ball $B_{\pp(l^2)}([w_{n_0}],\frac{1}{4})$ of $\pp(l^2)$. Since the distance is less that $1/4$ the sequence $([w_n])_{n\geqslant n_0}$ is in the same chart and we can suppose that the chart is $(\calu_0,h_0)$. Let denote $[w_n]=[1:\tilde w_n]$. 
Therefore $h_0(B_{\pp(l^2)}([w_{n_0}],\frac{1}{4}))$ is an open ball of $l^2$ that is complete and the sequence $h_0([w_n])$ converge to some $w$. Consider then $h_0^{-1}(w)=[1:w] \in \pp(l^2)$.
Then consider $n$ such that the path $\gamma: [0,1] \to \pp(l^2)$ defined by $\gamma(t)=[1 ; wt+(1-t)\tilde w_n]$ does not pass to $[1:0: \dots]$.  
Since 
\begin{equation}
	\begin{array}{lcl}
		\omega_{\gamma(t)}(\gamma'(t),\gamma'(t)) &=& \frac{\norm{\gamma(t)}^2 \norm{\gamma'(t)} - \langle \gamma(t),\gamma'(t)\rangle \langle\gamma'(t),\gamma(t)\rangle}{\norm{\gamma(t)}^4}  \\
		&=& \frac{1}{\norm{\gamma(t)}^4}\langle \ \norm{\gamma(t)}^2 \gamma'(t) - \langle \gamma'(t),\gamma(t) \rangle \gamma(t) \  , \ \gamma'(t) \ \rangle \\
		&\leqslant& \frac{\norm{\norm{\gamma(t)}^2 \gamma'(t) - \langle \gamma'(t),\gamma(t) \rangle \gamma(t)}}{\norm{\gamma(t)}^4} \norm{\gamma'(t)} \leqslant 2 \frac{\norm{\gamma'(t)}^2}{\norm{\gamma(t)}^2}
\end{array} \end{equation} 
one can write $\text{dist}([1:\tilde w_n],[1:w]) \leqslant C \norm{\tilde w_n - w}$ where $C>0$. Therefore the sequence $(w_n)_n$ converges to $[1:w]$ in $\pp(l^2)$ and it is complete.

\newsect[DIS-LOOP]{Distance on the loop space}
Let us consider $X$ a finite dimensional complex K\"ahler manifold with $h$ its Hermitian metric and $\nabla$ the unique connection compatible with $h$ (the Levi-Civita connection). In this context, a path $\gamma$ is a geodesic if it satisfies the geodesic equation 
\begin{equation}
	\nabla_{\dot{\gamma}(t)} \dot{\gamma}(t) = 0.
\end{equation}
Let $g \in LX$. For a vector field $\xi$ on $LX$ and a smooth function $\alpha$ on $LX$, there exist an induced vector field $\hat{\xi}$ on $X$  along $g$ and a smooth function $\hat \alpha$ on $X$ such that
\begin{equation}
	\forall s \in \ss^1, \, \hat{\xi}(g(s))=\xi_g(s) \text{ and } 
	\widehat{\alpha}(g(s))= \alpha(g).
\end{equation}

 Then we can define a connection $\tilde{\nabla}$ on the loop space $LX$ by the following. For two vector fields $\xi$ and $\nu$, we have $\tilde{\nabla}_\nu\xi : LX \to T(LX) $ such that
\begin{equation}
\forall  g \in LX, \,	\tilde{\nabla}_\nu\xi(g)= \Big(s \mapsto {\nabla}_{\hat\nu}\hat\xi(g(s)) \Big) \in T_g(LX) = W^{k,2}(\ss^1, g^*TX)
\end{equation}
\begin{prop}
	$\tilde{\nabla}$ is a connection on $LX$. 
\end{prop}
\proof Let $g \in LX$. 
For any smooth function $\alpha$ and vector fields $\xi$, $\eta$ on $LX$
\[
\forall s \in \ss^1, \, 
\left(\tilde{\nabla}_{\alpha \xi}\eta \right)_g(s) =  \nabla_{  \widehat{\alpha \xi}}\hat{\eta} (g(s)) = 
\hat\alpha(g(s))\nabla_{  \hat{\xi}}\hat{\eta} (g(s)) = \alpha(g) \left(\tilde{\nabla}_{\xi}\eta \right)_g(s)
\] 
It satisfies the Leibniz rule: Let $g \in LX$. For all $s \in \ss^1$
\begin{eqnarray*}
	\left(\tilde{\nabla}_{ \xi}(\alpha\eta) \right)_g(s) &=&  {\nabla}_{ \hat\xi}(\widehat{\alpha\eta})(g(s))=
	d_{g(s)}\hat\alpha \left[\hat\xi(g(s))\right]  \hat\eta(g(s))
	 + \hat{\alpha}(g(s)) \nabla_{  \hat{\xi}}\hat\eta (g(s)) \\
	 &=& d_g\alpha [\xi_g] \eta_g(s) + \alpha(g) \left(\tilde{\nabla}_\xi\eta \right)_g(s)
\end{eqnarray*}
\qed
\begin{prop}
	The connecion $\tilde{\nabla}$ satifies the two following properties. 
	\begin{enumerate}
		\item $\tilde{\nabla}$ is compatible with the metric $H$ defined in \eqqref(herm-m1).
		\item $\tilde{\nabla}$ is torsion free, i.e. $\nabla_\xi\eta - \nabla_\eta \xi -[\xi,\eta] =0$.
	\end{enumerate}
	
\end{prop}
\proof 
\begin{enumerate}
	\item For $g \in LX$ and $\xi, \eta, \nu$ vector fields on $LX$
	\[
	H_g(\tilde{\nabla}_\nu\xi, \eta) + H_g(\xi, \tilde{\nabla}_\nu\eta)
	= \int_{s\in \ss^1} h_{g(s)}\Big(\nabla_{\hat\nu}\hat\xi(g(s)) , \eta_g(s)\Big) + h_{g(s)}\Big(\xi_g(s), \nabla_{\hat\nu}\hat\eta(g(s))\Big)  \, ds
	\] 
	\[
	= \int_{s \in \ss^1} \hat\nu \Big(h_{g(s)}(\xi_g(s),\eta_g(s)) \Big)(g(s))ds 	= \int_{s \in \ss^1} \nu \Big(h_{g(s)}(\xi_g(s),\eta_g(s)) \Big)(g)ds 
	\]
	\begin{equation}
		= \nu \big(H(\xi,\eta) \big)(g) \quad \text{ by } \eqqref(diff-4)
	\end{equation}
	\item For $g \in LX$ and for all $s \in \ss^1$ 
	\[
	(\tilde\nabla_\xi\eta)_g(s) - (\tilde\nabla_\eta \xi)_g(s) -[\xi,\eta]_g(s)  = \nabla_{\hat{\xi}}\hat{\eta}(g(s)) - \nabla_{\hat{\eta}} \hat{\xi}(g(s)) -[\hat{\xi},\hat\eta](g(s)) = 0 
	\]
	Since $\nabla$ is the Levi-Civita connection on $X$. 
\end{enumerate}
\smallskip\qed
By the fundamental theorem of riemannian geometry, $\tilde{\nabla}$ is the Levi-Civita connection on $LX$. Thus we can state the following result. 
\begin{thm}
	Let $f$ and $g$ two loops on $LX$ and let $\gamma$ be a path joining them, i.e.
	\[
	\forall s \in \ss^1, \quad \gamma_s(0)= f(s) \text{ and } \gamma_s(1)=g(s).
	\]
	$\gamma$ is a geodesic on $LX$ if and only if for all $s \in \ss^1$, $\gamma_s$ is a geodesic on $X$ joining $f(s)$ and $g(s)$. 
\end{thm}
\proof 
The path $\gamma$ satisfies the geodesic equation : 
\[
\tilde{\nabla}_{\dot{\gamma}(t)}\dot{\gamma}(t)= 0 \quad \Longleftrightarrow \quad  \forall s \in \ss^1, \nabla_{\dot{\gamma}_s(t)}\dot{\gamma}_s(t)= 0 
\]
\qed

\ifx\undefined\bysame
\newcommand{\bysame}{\leavevmode\hbox to3em{\hrulefill}\,}
\fi

\def\entry#1#2#3#4\par{\bibitem[#1]{#1}
	{\textsc{#2 }}{\sl{#3} }#4\par\vskip2pt}{}{}{}

\end{document}